\documentclass[12pt]{article}
\usepackage[utf8]{inputenc}
\usepackage{amsmath}
\usepackage{amssymb}
\usepackage{amsthm}
\usepackage{fullpage}
\usepackage{graphicx}
\usepackage{hyperref}
\usepackage{color}
\usepackage[sc]{mathpazo}
\usepackage[T1]{fontenc}
\usepackage[margin=10pt,font=small]{caption}
\usepackage{enumitem}

\newtheorem{thm}{Theorem}[section]
\newtheorem{prop}[thm]{Proposition}
\newtheorem{lem}[thm]{Lemma}
\newtheorem{coro}[thm]{Corollary}
\newtheorem{conj}[thm]{Conjecture}
\theoremstyle{remark}
\newtheorem{rmk}{Remark}
\theoremstyle{definition}
\newtheorem{defn}[thm]{Definition}

\newcommand{\tdef}[1]{\textcolor{blue}{\emph{#1}}}
\newcommand{\horiz}{\begin{center}\rule{0.3\textwidth}{0.5pt}\end{center}}

\newcommand{\unary}{\operatorname{unary}}
\newcommand{\leaf}{\operatorname{leaf}}
\newcommand{\conn}{\mathcal{C}}
\newcommand{\redskl}{\mathcal{RS}}
\newcommand{\skl}{\mathcal{S}}
\newcommand{\sklf}{\operatorname{Sk}}
\newcommand{\syndiag}{\operatorname{Diag}}
\newcommand{\vtree}{\mathcal{V}}
\newcommand{\bipartite}{\mathcal{B}}
\newcommand{\planar}{\mathcal{M}}
\newcommand{\onecorner}{\mathcal{U}}

\newcommand{\excess}{\operatorname{ex}}
\newcommand{\uchain}{\operatorname{uc}}
\newcommand{\applv}{\operatorname{Appl}_{v}}
\newcommand{\appla}{\operatorname{Appl}_{a}}

\newcommand{\lnode}{\operatorname{lnode}}
\newcommand{\znode}{\operatorname{znode}}
\newcommand{\rlabel}{\operatorname{rlabel}}
\newcommand{\edge}{\operatorname{edge}}

\newcommand{\white}{\operatorname{white}}
\newcommand{\black}{\operatorname{black}}
\newcommand{\face}{\operatorname{face}}
\newcommand{\outdeg}{\operatorname{outdeg}}

\newcommand{\outvertex}{\operatorname{outv}}

\newcommand{\insertfig}[2]{\includegraphics[page=#1, width=#2\textwidth]{fig-ipe-ext.pdf}}

\author{Wenjie Fang \\ Univ Gustave Eiffel, CNRS, LIGM, F-77454 Marne-la-Vallée, France}

\title{Bijections between planar maps and planar linear normal $\lambda$-terms with connectivity condition}
\begin{document}

\maketitle

\abstract{The enumeration of linear $\lambda$-terms has attracted quite some attention recently, partly due to their link to combinatorial maps. Zeilberger and Giorgetti (2015) gave a recursive bijection between planar linear normal $\lambda$-terms and planar maps, which, when restricted to 2-connected $\lambda$-terms (\textit{i.e.}, without closed sub-terms), leads to bridgeless planar maps. Inspired by this restriction, Zeilberger and Reed (2019) conjectured that 3-connected planar linear normal $\lambda$-terms have the same counting formula as bipartite planar maps. In this article, we settle this conjecture by giving a direct bijection between these two families. Furthermore, using a similar approach, we give a direct bijection between planar linear normal $\lambda$-terms and planar maps, whose restriction to 2-connected $\lambda$-terms leads to loopless planar maps. This bijection seems different from that of Zeilberger and Giorgetti, even after taking the map dual. We also explore enumerative consequences of our bijections.}

\horiz

As a well-known Turing-complete computational model, $\lambda$-calculus has been well-studied in logic and related fields. However, their enumeration did not attract attention in combinatorics until relatively recently. One of the most well-studied family is \textit{linear $\lambda$-terms} (also called \textit{BCI terms} from combinatorial logic), due to their connections with combinatorial maps, which are graph embeddings on surfaces. Such connections were pioneered by Bodini, Gardy and Jacquot in \cite{lambda-cubic}, where a simple bijection between linear $\lambda$-terms and cubic maps was given. Due to its simplicity, this bijection transfers many interesting statistics from linear $\lambda$-terms to cubic maps, and can be specialized to interesting sub-families. This approach has also led to the study of asymptotic properties and statistics distribution of related $\lambda$-terms (see, \textit{e.g.}, \cite{lambda-cubic, generalized-BCI, linear-stat}), and further connections to other objects \cite{chord-diagram-bridgeless}.

Independently, Zeilberger and Giorgetti noticed in \cite{planar-lambda} a new connection between a sub-family of linear $\lambda$-terms called \emph{planar linear normal $\lambda$-terms} and planar maps. In fact, they gave a bijection between the two families, which is not a simple restriction of the one in \cite{lambda-cubic} for general linear $\lambda$-terms. When restricted to $\lambda$-terms that are also \emph{unitless}, \textit{i.e.}, without closed sub-terms, this bijection leads to bridgeless planar maps. The unitless condition here is equivalent to the 2-connectedness of the syntactic diagram of the $\lambda$-term. Such connections lead naturally to the consideration of higher connectivity conditions on planar linear $\lambda$-terms. In a talk at CLA 2019 \cite{3-conn-conj}, Zeilberger and Reed considered possible interpretation of connectivity in logic, using the formalism in \cite{imploid}. Based on computer experiments, they also proposed the following conjecture.

\begin{conj}[\cite{3-conn-conj}] \label{conj:3-conn}
  The number of 3-connected planar linear normal $\lambda$-terms with $n+2$ variables is
  \[
    \frac{2^n}{(n+1)(n+2)} \binom{2n+1}{n},
  \]
  which is also the number of bipartite planar maps with $n$ edges \cite{Tutte-census}.
\end{conj}

In a talk at CLA 2020 \cite{grygiel-yu}, Grygiel and Yu tried to relate these 3-connected $\lambda$-terms to $\beta(0,1)$-trees, which are in bijection with bicubic planar maps \cite{desc-tree} and thus with bipartite planar maps \cite{Tutte-census}, and succeeded in a few special cases. The key of their work is a characterization of skeletons of such $\lambda$-terms (Proposition~\ref{prop:chara-lambda}).

In this article, based on the characterization given by Grygiel and Yu in \cite{grygiel-yu}, we prove Conjecture~\ref{conj:3-conn} by giving a bijection between 3-connected planar linear normal $\lambda$-terms and bipartite planar maps, passing through the so-called \emph{degree trees} defined in \cite{fang-bipartite}. We have the following theorem.

\begin{thm} \label{thm:main-3-conn}
  For $n \geq 2$, there is a bijection between 3-connected planar linear normal $\lambda$-terms with $n$ variables and bipartite planar maps with $n-2$ edges.
\end{thm}

Using this new bijection, we study the refined enumeration of such $\lambda$-terms under several statistics, and also the asymptotic distribution of some statistics. All is done by translating known results on bipartite planar maps.

With a similar approach, we also give a new bijection between planar linear normal $\lambda$-terms and planar maps, which is direct and also has a nice restriction to 2-connected $\lambda$-terms. As a direct bijection, it may be a better tool to explore the relation between planar maps and related $\lambda$-terms.

\begin{thm} \label{thm:main-conn}
  There is a direct bijection between planar linear normal $\lambda$-terms with $n$ variables and planar maps with $n-1$ edges. Furthermore, it sends those $\lambda$-terms that are also 2-connected to loopless planar maps.
\end{thm}

For this bijection, we define a family of node-labeled trees called \emph{v-trees}, which can be seen as description trees (see \cite{desc-tree}) of planar maps under a seemingly new recursive decomposition called \emph{one-corner decomposition}. We then give direct bijections from v-trees to both planar maps and planar linear normal $\lambda$-terms, thus linking the two families. Our bijection is different from that in \cite{planar-lambda}, even after taking the map dual, which sends loopless planar maps to bridgeless planar maps.

The rest of the article is organized as follows. In Section~\ref{sec:prelim}, we introduce the families of $\lambda$-terms and planar maps that will be our subject of study. Then we characterize in Section~\ref{sec:chara} the skeletons of planar linear normal $\lambda$-terms with different connectivity conditions. Using these characterizations, in Section~\ref{sec:bij-3} we give a bijection from 3-connected planar linear normal $\lambda$-terms to degree trees, which are in bijection with bipartite maps. We also study some refined enumerations and the asymptotic distribution of some statistics of the related $\lambda$-terms. Using the same approach, in Section~\ref{sec:bij-2}, we define v-trees and relate them bijectively to planar linear normal $\lambda$-terms. Then finally, in Section~\ref{sec:map-decomp}, we propose the one-corner decomposition of planar maps, and show that v-trees describe this recursive decomposition. We also give a direct bijection between planar maps and v-trees, completing our new bijection between planar maps and planar linear normal $\lambda$-terms.

\paragraph{Acknowledgment} We thank Katarzyna Grygiel and Guan-Ru Yu for discovering Proposition~\ref{prop:chara-lambda} and for their inspiring discussions. We also thank \'Eric Fusy for showing how the decomposition of bicubic planar maps in \cite{desc-tree} can be translated to bipartite maps, inspiring the one-corner decomposition. We are grateful to the organizers of the workshop series of Computational Logic and Applications (CLA), in which Grygiel and Yu got to know Conjecture~\ref{conj:3-conn} from Zeilberger and Reed, and in which the author got to know Proposition~\ref{prop:chara-lambda} from Grygiel and Yu. We are also grateful for the useful comments from Noam Zeilberger. The author is partially supported by ANR LambdaComb (ANR-21-CE48-0017).

\section{Preliminaries} \label{sec:prelim}

In mathematical logic, \tdef{$\lambda$-terms} can be defined recursively using three constructions:
\begin{itemize}
\item Atoms, which are variables $x, y, \ldots$;
\item Application $t \; u$ with two $\lambda$-terms $t$ and $u$;
\item Abstraction $\lambda x.t$ with $t$ a $\lambda$-term and $x$ a variable. We say that every atom with the variable $x$ in $t$ is \tdef{bound} by the abstraction, if not yet bound.
\end{itemize}
In the following, we consider $\lambda$-terms up to \tdef{$\alpha$-renaming}, \textit{i.e.}, changing the name of a variable in an abstraction consistently throughout all atoms it binds in the $\lambda$-term. A $\lambda$-term $t$ is \tdef{closed} if all atoms are bound, and it is \tdef{linear} if it is closed and each abstraction binds exactly one atom.

Computation on $\lambda$-terms is defined by the \tdef{$\beta$-reduction}: if in a $\lambda$-term $t$, there is a sub-term of the form $(\lambda x. u) v$, then we may perform a $\beta$-reduction on it for $t$ to obtain another $\lambda$-term $t'$. This is done by replacing $(\lambda x. u) v$ with $u[x \leftarrow v]$, the sub-term $u$ with all the atoms of $x$ replaced by a copy of the sub-term $v$. A $\lambda$-term $t$ is \tdef{$\beta$-normal} (or simply \tdef{normal}) if no $\beta$-reduction is possible for $t$.

In the context of enumeration, we usually consider a graphical presentation of $\lambda$-terms based on its tree structure. Given a $\lambda$-term $t$, its \tdef{skeleton} $\sklf(t)$ is a plane unary-binary tree defined recursively as follows:
\begin{itemize}
\item If $t$ is an atom, then $\sklf(t)$ is a leaf;
\item If $t$ is an application $u \; v$, then $\sklf(t)$ is rooted at a binary node, with $\sklf(u)$ its left sub-tree, and $\sklf(v)$ its right sub-tree;
\item If $t$ is an abstraction $\lambda x.u$, then $\sklf(t)$ is rooted at a unary node, with $\sklf(u)$ its only sub-tree.
\end{itemize}
$\sklf(t)$ is simply the syntactic tree of $t$, and we identify atoms and abstractions in $t$ with their corresponding leaves and unary nodes in $\sklf(t)$. Many properties of $\lambda$-terms can be read off directly from their skeletons. For instance, a $\lambda$-term $t$ is normal if no binary node has a unary node as its left child. The \tdef{size} of a skeleton is the number of its leaves, which is also the number of atoms in the original $\lambda$-term.

\begin{figure}
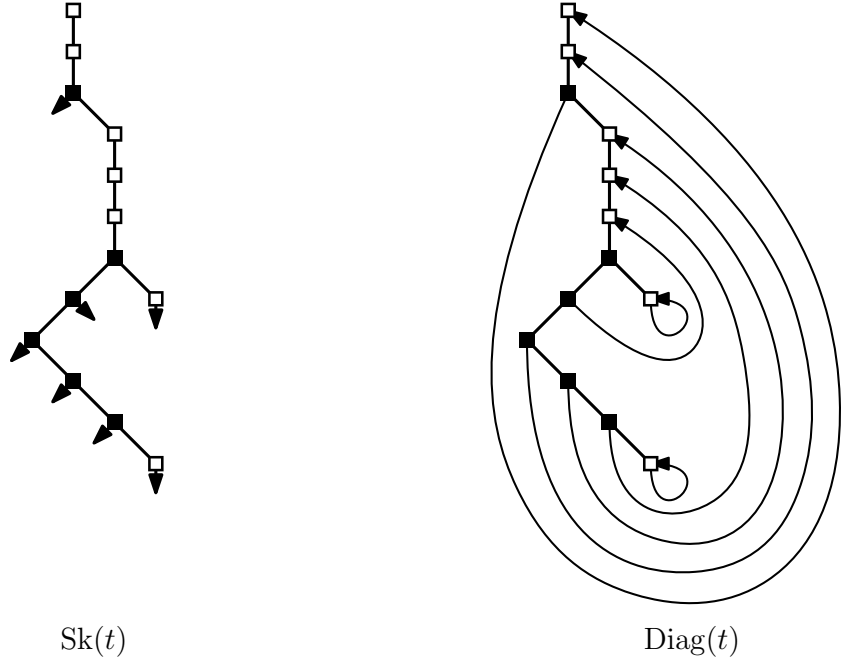

  \centering
  \insertfig{6}{0.67}
  \caption{Example of a normal planar $\lambda$-term, its skeleton and syntactic diagram.}
  \label{fig:lambda-example}
\end{figure}

The skeleton $\sklf(t)$ of a $\lambda$-term $t$ still misses the binding relations between atoms and abstractions. The \tdef{syntactic diagram} (or simply \tdef{diagram}) of a $\lambda$-term $t$, denoted by $\syndiag(t)$, is obtained from its skeleton $\sklf(t)$ by replacing each leaf by an edge from its parent to the unary node of abstraction that binds the atom of the leaf. It is clear that two different $\lambda$-terms (up to $\alpha$-renaming) never share the same diagram. Given a leaf of $\sklf(t)$, we draw the corresponding edges starting from the parent of the leaf, traveling in counter-clockwise direction and entering the unary node from the right. If such a drawing can be done without intersection, then we say that the $\lambda$-term $t$ is \tdef{planar}. Figure~\ref{fig:lambda-example} provides an example of a planar $\lambda$-term, with its skeleton and syntactic diagram. We note that the $\lambda$-term in the example is also normal.

As diagrams are graphs, we may transfer notions in graph theory to $\lambda$-terms. If the diagram of a $\lambda$-term $t$ is \tdef{2-edge-connected} (or simply \tdef{2-connected}), \textit{i.e.}, the removal of any edge does not disconnect the diagram, then we say that $t$ is \tdef{2-connected}. Similarly, a $\lambda$-term $t$ is \tdef{3-edge-connected} (or simply \tdef{3-connected}) if the removal of any two edges in $\syndiag(t)$ does not disconnect $\syndiag(t)$, except eventually the two edges adjacent to the root of $\sklf(t)$ when $t$ starts with an abstraction binding one atom.

In this article, we consider planar linear $\lambda$-terms, in which each atom is bound by a different ancestral abstraction. For planarity, given the skeleton of such a $\lambda$-terms, we take a clockwise contour walk starting from the root of the skeleton. In the walk, we read off the unary nodes and leaves by their first visit, and by planarity, we must obtain a well-parenthesized word with unary nodes acting as opening parentheses and leaves as closing ones, with each abstraction paired up with the atom it binds. Therefore, we may identify planar linear $\lambda$-terms with their skeletons, as the way atoms are bound is given by planarity. By abuse of notation, given a skeleton $S$ of some planar linear $\lambda$-terms $t$, we also denote by $\syndiag(S)$ the diagram $\syndiag(t)$ of $t$.

\begin{rmk}
  We note that, in the definition of planar $\lambda$-term, there are two choices of the drawing, either going clockwise or counter-clockwise from the parent of the leaf to the corresponding unary node. The two choices define different families of planar $\lambda$-terms, although they are all in bijection with the same set of unary-binary trees. As we will see, the difference of the two families does not matter to connected and 2-connected planar linear $\lambda$-terms, but it matters to 3-connected terms. We choose the current definition (also called RL-planarity), as it is the one used in Conjecture~\ref{conj:3-conn} \cite{3-conn-conj}. This convention is more natural in the sense that it is stable under $\beta$-reduction, which the other is not. For a discussion on the two conventions, see \cite{planar-lambda}.
\end{rmk}

In the following, we relate previously defined $\lambda$-terms to planar maps. A \tdef{planar map} is an embedding of a graph on the plane such that edges only intersect at vertices, defined up to orientation-preserving continuous deformations. The connected components of the complement of the embedding are called the \tdef{faces} of the planar map. The face that extends to infinity is called the \tdef{outer face}. Furthermore, we only consider \tdef{rooted} planar maps here, which means we mark a special corner on the outer face, called the \tdef{root corner}. The edge next to the root corner in the clockwise order is called the \tdef{root edge}, and the vertex of the corner is called the \tdef{root vertex}. We also consider the \tdef{empty map}, which consists of a single vertex. A planar map is \tdef{bipartite} if we can color its vertices with black and white such that each edge links a black and a white vertex. The root vertex is colored black by convention.

\section{Characterization of planar normal $\lambda$-terms with connectivity conditions} \label{sec:chara}

As mentioned in the previous section, we may classify planar linear $\lambda$-terms by the connectivity of their diagram. We denote by $\conn^{(1)}_n$ the set of planar linear normal $\lambda$-terms of size $n$. Similarly, we denote by $\conn^{(2)}_n$ (resp. $\conn^{(3)}_n$) the terms that are 2-connected (resp. 3-connected) in $\conn^{(1)}_n$. We take $\conn^{(1)} = \bigcup_{n \geq 1} \conn^{(1)}_n$, and we define $\conn^{(2)}$ and $\conn^{(3)}$ similarly. It is clear that $\conn^{(3)} \subset \conn^{(2)} \subset \conn^{(1)}$.

Due to the correspondence between planar linear $\lambda$-terms and their skeletons, in the following we will focus on skeletons instead of $\lambda$-terms. We denote by $\skl^{(1)}_n$ (resp. $\skl^{(2)}_n$ and $\skl^{(3)}_n$) the set of skeletons of terms in $\conn^{(1)}_n$ (resp. $\conn^{(2)}_n$ and $\conn^{(3)}_n$). We also take $\skl^{(1)} = \bigcup_{n \geq 1} \skl^{(1)}_n$, and we define $\skl^{(2)}$ and $\skl^{(3)}$ similarly. It is clear that $\skl^{(3)} \subset \skl^{(2)} \subset \skl^{(1)}$.

We start by characterizing skeletons of planar linear $\lambda$-terms. Given a unary-binary tree $S$ and a node $u$ in $S$, we denote by $S_u$ the sub-tree induced by $u$, by $\leaf(S)$ the number of leaves in $S$, and by $\unary(S)$ the number of unary nodes in $S$. We also define the \tdef{excess} of $S$, denoted by $\excess(S)$, to be $\leaf(S) - \unary(S)$.

\begin{prop} \label{prop:chara-conn-lambda}
  A unary-binary tree $S$ is in $\skl^{(1)}$ if and only if
  \begin{itemize}
  \item (\textbf{Linearity}) $\excess(S) = 0$;
  \item (\textbf{Normality}) No binary node in $S$ has a unary node as its left child;
  \item (\textbf{Connectedness}) For every binary node and leaf $u$ in $S$, we have $\excess(S_u)$ greater than or equal to the number of consecutive unary nodes above $u$.
  \end{itemize}
  Furthermore, a unary-binary tree $S \in \skl^{(1)}$ is also in $\skl^{(2)}$ if and only if
  \begin{itemize}
  \item (\textbf{2-connectedness}) For every binary node and leaf $u$ in $S$, we have $\excess(S_u)$ strictly larger than the number of consecutive unary nodes above $u$.
  \end{itemize}
\end{prop}
\begin{proof}
  We first consider the $\skl^{(1)}$ case. The ``only if'' part comes directly from the definition of terms in $\conn^{(1)}$: (\textbf{Linearity}) from the linearity of the term, (\textbf{Normality}) from $\beta$-normality, and (\textbf{Connectedness}) from the fact that every abstraction binds a variable in the sub-term. For the ``if'' part, given $S$ satisfying the conditions, we can pair unary nodes with leaves in a unique way to form a linear and planar term, and such a process is guaranteed to succeed by (\textbf{Connectedness}).

  For the $\skl^{(2)}$ case, we first observe that (\textbf{2-connectedness}) implies (\textbf{Connectedness}). Then, the sub-term rooted at the top of the chain of unary nodes directly above $u$ is closed if and only if there are exactly $\excess(S_u)$ consecutive unary nodes above $u$, meaning that (\textbf{2-connectedness}) is equivalent to the absence of closed sub-terms. We now show that the absence of closed sub-terms of a $\lambda$-term $t$ in $\conn^{(1)}$ is equivalent to $t \in \conn^{(2)}$, which is in turn equivalent to $\sklf(t) \in \skl^{(2)}$. Consider the diagram $\syndiag(t)$ of $t$, which is built on its skeleton $\sklf(t)$. Removing any edge in $\syndiag(t)$ coming from a leaf of $\sklf(t)$ will not disconnect $\syndiag(t)$, as edges in $\sklf(t)$ still connects all nodes of $\syndiag(t)$. If removing an edge $e$ of $\sklf(t)$ disconnects $\syndiag(t)$, then the sub-tree of $\sklf(t)$ induced by $e$ must correspond to a closed sub-term, as no leaf is transformed into an edge linking to any unary node in $\sklf(t)$ above $e$. Conversely, a closed sub-term in $t$  corresponds to a sub-tree $S_u$ of $\sklf(t)$ rooted at a node $u$, where all leaves of $S_u$ are transformed into edges in $\syndiag(t)$ that link nodes in $S_u$. Therefore, cutting the edge leading to $u$ disconnects $\syndiag(t)$. We thus have the claimed equivalence.
\end{proof}

From the proof above, we see that a term $t \in \conn^{(1)}$ is also in $\conn^{(2)}$ if it is unitless, \textit{i.e.}, has no closed sub-terms.

For skeletons of $\lambda$-terms in $\conn^{(3)}$, we start with an observation from \cite{grygiel-yu}.

\begin{prop}[See \cite{grygiel-yu}] \label{prop:lambda-root}
  Let $S \in \skl^{(3)}$ and $u$ its first binary node. Then the left child of $u$ is a leaf.
\end{prop}
\begin{proof}
  Let $v$ be the left child of $u$. If $v$ is not a leaf, as $S \in \skl^{(3)} \subset \skl^{(2)}$, it induces a sub-term with at least one free variable. By planarity, the leafs of $S_v$ corresponding to free variables in the sub-term are connected to an initial segment of the consecutive unary nodes above $u$. Cutting $\{u, v\}$ and the edge immediately below the said initial segment disconnects $\syndiag(S)$, contradicting 3-connectedness.
\end{proof}

By Proposition~\ref{prop:lambda-root}, we define the \tdef{reduced skeleton} of $t \in \conn^{(3)}$ to be the right sub-tree of the first binary node in $S = \skl(t)$. By Proposition~\ref{prop:lambda-root} and (\textbf{Linearity}) on $S$, we can reconstruct $S$ from the reduced skeleton of $t$. We denote by $\redskl_n$ the set of reduced skeletons of terms in $\conn^{(3)}_n$, and $\redskl = \cup_{n \geq 3} \redskl_n$. The following characterization of reduced skeletons was first discovered and stated without proof in \cite{grygiel-yu}.

\begin{prop}[See \cite{grygiel-yu}] \label{prop:chara-lambda}
  A unary-binary tree $S$ is in $\redskl$ if and only if
  \begin{itemize}
  \item (\textbf{Normality}) No binary node in $S$ has a unary node as its left child;
  \item (\textbf{3-connectedness}) For every binary node $u$ in $S$, let $v$ be its right child, then $\excess(S_v)$ is strictly larger than the number of consecutive unary nodes above $u$.
  \end{itemize}
\end{prop}
\begin{proof}
  The condition (\textbf{Normality}) comes directly from the definition of normal $\lambda$-terms. If (\textbf{3-connectedness}) fails for some binary node $u$ with $v$ its right child, then by planarity, the leafs in $S_v$ are connected to some consecutive unary nodes directly above $u$. Thus, by cutting the edge from $u$ to its \emph{left child} and the edge immediately above the said consecutive unary nodes, we disconnect the diagram of the corresponding $\lambda$-term, meaning that $S \notin \redskl$.

  Now, conversely, suppose that the two conditions holds. We first show that (\textbf{3-connectedness}) implies (\textbf{2-connectedness}) in Proposition~\ref{prop:chara-conn-lambda}. For a binary node $u$ in $S$, let $u_0 = u, u_1, \ldots$ be the nodes such that $u_{i+1}$ is the left child (or the only child) of $u_i$. It ends at some node $u_k$, which is a leaf. By (\textbf{Normality}), all $u_i$ except $u_k$ are binary nodes. Let $v_i$ be the right child of $u_i$ for $0 \leq i < k$. Suppose that there are $\ell$ unary nodes above $u$. Then we have
  \[
    \excess(S_u) = \leaf(S_u) - \unary(S_u) - \ell = 1 + \sum_{i = 0}^{k - 1} \leaf(S_{v_i}) - \ell - \sum_{i = 0}^{k-1} \unary(S_{v_i}) > 1.
  \]
  The last inequality comes from applying (\textbf{3-connectedness}) on all $u_i$ and observing that there is no unary node directly above any $u_i$ for $i \geq 1$. Thus, (\textbf{2-connectedness}) holds for all binary nodes in $S$. Now for leaves in $S$, if it is at the end of a chain of unary nodes, then by (\textbf{Normality}), the first node of the chain must be the right child of a binary node, and then (\textbf{3-connectedness}) implies that there is at most 0 unary nodes in the chain, which is impossible. Therefore, every leaf in $S$ is a child of a binary node, in which case (\textbf{2-connectedness}) holds. 

  Applying (\textbf{2-connectedness}) to the first binary node of $S$, we have $\excess(S) > 0$. Let $S^+$ be the unary-binary tree that starts with $\excess(S) + 1$ unary nodes, then a binary node $u_+$ with a leaf as left child and $S$ as right sub-tree. By construction, $S^+$ satisfies (\textbf{Linearity}) in Proposition~\ref{prop:chara-conn-lambda}. We also check that (\textbf{2-connectedness}) holds for $u_+$ and the new added leaf. By Proposition~\ref{prop:chara-conn-lambda}, we thus conclude that $S^+ \in \skl^{(2)}$, meaning that $\syndiag(S^+)$ is 2-connected.

  Now we take an arbitrary pair of edges $e_1, e_2$ of $\syndiag(S^+)$, not both adjacent to the root, and we show that their removal does not disconnect $\syndiag(S^+)$. We proceed by case analysis on the origin of the two edges.
  \begin{itemize}
  \item \textbf{Both $e_1, e_2$ come from leaves in $S^+$.} \\
    In this case, edges in $S^+$ hold $\syndiag(S^+)$ connected.
  \item \textbf{One of $e_1, e_2$ is an edge of $S^+$ and the other comes from a leaf.} \\
    We suppose that $e_1$ is an edge of $S^+$. Let $u$ be the first binary node in the sub-tree of $S^+$ induced by cutting $e_1$. If $u$ is the newly added binary node $u_+$ in $S^+$, then as all nodes above it are unary nodes, which are turned into trivalent nodes in $\syndiag(S^+)$ except the root. The only possibility that $e_1, e_2$ disconnect $\syndiag(S^+)$ in this sub-case is that both $e_1, e_2$ are adjacent to the root, which is discarded in the definition. Otherwise, $u$ is also in $S$. Let $v$ be the right child of $u$. Suppose that there are $\ell$ consecutive unary nodes directly above $u$. By (\textbf{3-connectedness}) on $u$, we have $\excess(S_v) \geq \ell + 1$, meaning that there is at least one leaf in $S_v$ that is transformed into an edge $e_v$ in $\syndiag(S^+)$ linking a node in $S_v$ to somewhere above the consecutive unary nodes above $u$, thus also above $e_1$. By (\textbf{Normality}), there is no unary node on the path from $u$ to the leftmost leaf of $S_u$. By planarity, given $e_v$, the edge in $\syndiag(S^+)$ corresponding to the leftmost leaf of $S_u$ links the parent of the leaf to somewhere above $e_1$. There are thus two edges in $\syndiag(S^+)$ from leaves in $S^+$ that link the part of $S^+$ below $e_1$ to the part above it. It is therefore impossible to disconnect $\syndiag(S^+)$ by further removing $e_2$ that comes from a leaf.
  \item \textbf{Both $e_1, e_2$ are edges of $S^+$.} \\
    If the two sub-trees $S_1, S_2$ of $S^+$ induced by cutting $e_1$ and $e_2$ respectively are not contained one in the other, then by (\textbf{2-connectedness}), in both $S_1$ and $S_2$, there is an edge in $\syndiag(S^+)$ that links some node to a unary node in the part of $S^+$ above both $e_1$ and $e_2$, meaning that the removal of $e_1, e_2$ cannot disconnect $\syndiag(S^+)$. We thus suppose that $e_1$ is in the sub-tree of $S^+$ induced by cutting $e_2$. Let $V_\uparrow$ be the set of nodes of $S^+$ above $e_2$, and $V_\downarrow$ be the set of nodes in the sub-tree of $S^+$ induced by cutting $e_1$. We denote by $V_{\mathrm{mid}}$ the set of nodes in $S^+$ that are not in either $V_\uparrow$ or $V_\downarrow$. This partition of nodes of. $S^+$ into $V_\uparrow, V_{\mathrm{mid}}, V_\downarrow$ is the one induced by the removal of $e_1, e_2$ in $S^+$. For $\syndiag(S^+)$, the removal of $e_1, e_2$ only disconnects it into two parts, if ever, as $S^+ \in \skl^{(2)}$. Furthermore, in this case, $V_\uparrow$ is connected to $V_\downarrow$ after the removal of $e_1, e_2$, otherwise the removal of one of $e_1$ and $e_2$ would already disconnect $\syndiag(S^+)$. We now show that it is always possible to find a leaf of $S^+$ that is transformed into an edge linking some node in $V_{\mathrm{mid}}$ to some node in $V_\uparrow \cup V_\downarrow$. Let $u$ be the upper end of $e_1$. We have $u \in V_{\mathrm{mid}}$. There are now three cases:
    \begin{itemize}
    \item \textbf{$u$ is a unary node.} \\
      In this case, as $S^+$ satisfies (\textbf{Linearity}), there is a leaf below $e_1$ (thus with its parent in $V_\downarrow$) corresponding to $u$, and in $\syndiag(S^+)$ it becomes an edge linking its parent in $V_\downarrow$ to $u \in V_{\mathrm{mid}}$.
    \item \textbf{$u$ is a binary node, and $e_1$ is the edge linking $u$ to its left child.} \\
      In this case, let $v$ be the right child of $u$, then all nodes in $S_v$ are in $V_{\mathrm{mid}}$. By (\textbf{3-connectedness}) on $u$, there is at least one leaf in $S_v$ that is transformed in $\syndiag(S^+)$ into an edge linking its parent, which is in $V_{\mathrm{mid}}$, to a unary node $w$ above the consecutive chain of unary nodes above $u$. If $w \in V_\uparrow$, then we have found the wanted leaf. Otherwise, we have $w \in V_{\mathrm{mid}}$. Let $u'$ be the first binary node below $w$. We have $u' \neq u$ and $u' \in V_{\mathrm{mid}}$. By (\textbf{2-connectedness}) on $u'$, the leftmost leaf in $S_{u'}$, whose parent is also in $V_{\mathrm{mid}}$, is transformed in $\syndiag(S^+)$ into an edge linking its parent to a unary node $w'$ above the consecutive chain of unary nodes above $u'$. We may thus repeat the whole process until reaching a node in $V_\uparrow$. The process will terminate, as we have only finitely many nodes in $S^+$. The leaf leading to the termination is the one we want.
    \item \textbf{$u$ is a binary node, and $e_1$ is the edge linking $u$ to its right child.} \\
      This case is settled in the same way as the previous one, except that we use (\textbf{2-connectedness}) on $u$ and we start with the leftmost leaf in $S_u$.
    \end{itemize}
    We thus show that, in any case, we can find such a leaf in $S^+$ that will be transformed in $\syndiag(S^+)$ into an edge linking some node in $V_{\mathrm{mid}}$ to some other node in $V_\uparrow \cup V_\downarrow$, meaning that the removal of $e_1, e_2$ cannot disconnect $\syndiag(S^+)$.
  \end{itemize}

  We conclude that, given (\textbf{Normality}) and (\textbf{3-connectedness}), the removal of any pair of edges $e_1, e_2$ cannot disconnect $\syndiag(S^+)$, meaning that $S^+ \in \skl^{(3)}$ and $S \in \redskl$. We thus have the claimed equivalence.
\end{proof}

\section{Bijection from $\lambda$-terms in $\conn^{(3)}$ to degree trees and refined enumeration} \label{sec:bij-3}

A \tdef{degree tree}, as defined in \cite{fang-bipartite}, is a pair $(T, \ell)$ where $T$ is a plane tree and $\ell$ a labeling on the nodes of $T$ such that
\begin{itemize}
\item If $u$ is a leaf, then $\ell(u) = 0$;
\item For $u$ with children $v_1, v_2, \ldots, v_k$ from left to right, let $s(u) = k + \sum_{i=1}^{k} \ell(v_i)$, then $s(u) - \ell(v_1) \leq \ell(u) \leq s(u)$.
\end{itemize}
The \tdef{size} of a degree tree $(T, \ell)$ is the number of edges in $T$, also denoted by $|T|$. We denote by $\mathcal{T}_n$ the set of degree trees of size $n$. Given a degree tree $(T, \ell)$, we define its edge labeling $\ell_\Lambda$ as follows: for each internal node $u$, we put a label $s(u) - \ell(u)$ on its leftmost descending edge, and $0$ on all other edges. It is clear that we can recover $(T, \ell)$ from $T$ and $\ell_\Lambda$. Later, we will use both $\ell$ and $\ell_\Lambda$. We use the following result from \cite{fang-bipartite} to convert $\ell_\Lambda$ to $\ell$.

\begin{lem}[See Lemma~1.1 in \cite{fang-bipartite}] \label{lem:degree-tree}
  For $T$ a plane tree and $\ell_\Lambda$ an edge labeling on $T$, the associated node labeling $\ell$ is given by $\ell(v) = |T_v| - \sum_{e \in T_v} \ell_\Lambda(e)$.
\end{lem}

We now define a bijection $\varphi$ from reduced skeletons to degree trees.

\begin{figure}
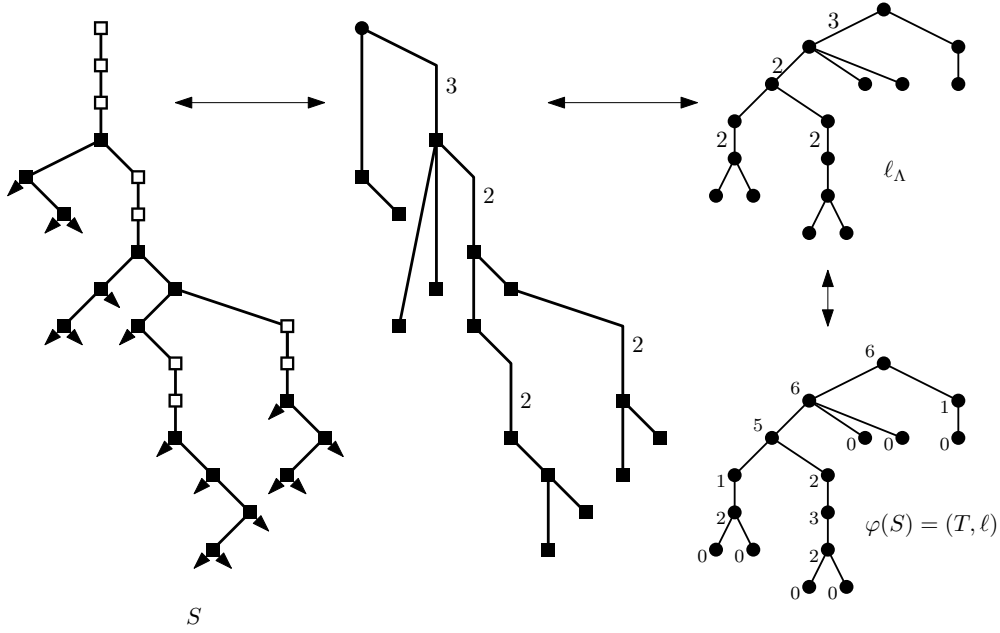

  \centering
  \insertfig{1}{0.8}
  \caption{Example of the bijection $\varphi$ from reduced skeletons in $\redskl$ to degree trees (zeros in $\ell_\Lambda$ are omitted).}
  \label{fig:bij-example}
\end{figure}

\begin{defn} \label{defn:redskl-degree-tree}
Given $S \in \redskl$, we perform the following:
\begin{enumerate}
\item Remove all leaves from $S$, then smooth out all unary nodes and put the number of smoothed-out nodes at each edge. For unary nodes above the first binary node, before the smoothing out, we add a new node as the new root, with the old root its right child. We obtain a binary tree with edge labels $(S', \ell_\Lambda')$.
\item Perform the classical bijection (\emph{switching left and right}) from binary tree $S'$ to plane tree $T'$: given a node $u$, its right child in $S'$ turns into its leftmost child in $T'$, and its left child in $S'$ into its sibling immediately to the right. The edge labels of $\ell_\Lambda'$, which are all on right edges in $S'$, become an edge labeling $\ell_\Lambda$ with only non-zero values on the leftmost descending edges of nodes.
\end{enumerate}
We define $\varphi(S)$ to be $(T', \ell)$, where $\ell$ is the node labeling corresponding to $\ell_\Lambda$ constructed above. See Figure~\ref{fig:bij-example} for an example of the bijection $\varphi$.
\end{defn}

For simplicity, given $S \in \redskl$ and $(T, \ell) = \varphi(S)$, we identify non-root nodes in $T$ with binary nodes in $S$. For the root of $T$, we may identify it with the first binary node in the skeleton corresponding to $S$. Now we show some properties of $\varphi(S)$.

\begin{prop} \label{prop:bij-stat}
  Let $S \in \redskl_n$ and $(T, \ell) = \varphi(S)$. Then we have
  \begin{enumerate}
  \item $T$ is of size $n-2$;
  \item For each node $u$, let $v$ be its right child in $S$, then $\unary(S_v) = \sum_{e \in T_u} \ell_\Lambda(e)$;
  \item For each node $u$, let $v$ be its right child in $S$, then $\leaf(S_v) = |T_u| + 1$.
  \end{enumerate}
\end{prop}
\begin{proof}
  For the first property, we note that $S \in \redskl_n$ has $n-1$ leaves, thus $n-2$ binary nodes, each leading to an edge in $T$. The second one comes directly from the construction. For the third, we see from the construction that nodes in $S_v$ are exactly nodes in $T_u \setminus \{u\}$. Moreover, $\leaf(S_v)$ is the number of binary nodes in $S_v$ plus one. Thus, $\leaf(S_v) = |T_u| + 1$, as edges in $T_u$ are in bijection with nodes in $T_u \setminus \{u\}$.
\end{proof}

Now we show that $\varphi(S)$ is always a degree tree.

\begin{prop} \label{prop:bij-redskl-degree-tree}
  A unary-binary tree $S$ is in $\redskl_n$ if and only if $\varphi(S)$ is in $\mathcal{T}_{n-2}$.
\end{prop}
\begin{proof}
  The size parameter of $\varphi(S)$ is given by Proposition~\ref{prop:bij-stat}. We only need to show that the conditions for $S$ to be in $\redskl$, given in Proposition~\ref{prop:chara-lambda}, are equivalent to those for $\varphi(S)$ to be a degree tree.

  Let $(T, \ell) = \varphi(S)$ and $\ell_\Lambda$ its edge labeling. From the construction, we see that the (\textbf{Normality}) for $S$ is equivalent to say that $\ell_\Lambda$ only has non-zero labels on leftmost descending edges for each node. Then, by Proposition~\ref{prop:bij-stat} and Definition~\ref{defn:redskl-degree-tree}, (\textbf{3-connectedness}) for $S$ is equivalent to say that, for each node $w$ in $T$, let $u$ be its leftmost child in $T$ and $e_*$ the edge between $w$ and $u$, and suppose that $v$ is the right child of $u$ in $S$, we have
  \[
    \excess(S_v) = \leaf(S_v) - \unary(S_v) = |T_u| + 1 - \sum_{e \in T_u} \ell_\Lambda(e) > \ell_\Lambda(e_*).
  \]
  By Lemma~\ref{lem:degree-tree}, this is equivalent to
  \[
    \ell_\Lambda(e_*) \leq |T_u| - \sum_{e \in T_u} \ell_\Lambda(e) = \ell(u),
  \]
  which is exactly the definition of degree trees, as $\ell_\Lambda(e_*) \geq 0$ by construction.
\end{proof}

Given the bijection from $\redskl$ to degree trees, which are also in bijection with bipartite planar maps, we can transfer enumeration results from maps to $\lambda$-terms via degree trees. Some statistics are also transferred alongside.

We first define several statistics on reduced skeletons. Given $S \in \redskl_n$, we denote by $\applv(S)$ the number of binary nodes whose right child is a leaf, $\appla(S)$ the number of binary nodes whose right child is also a binary node, and $\uchain_k(S)$ the number of maximal chains of unary nodes of length $k$. In the language of $\lambda$-terms, the excess $\excess(S)$ is the number of abstractions at the beginning of the original term minus one (from the deleted leaf in Proposition~\ref{prop:lambda-root}), $\applv(S)$ the number of applications of terms to a variable, $\appla(S)$ the number of applications of terms to another application, and $\uchain_k(S)$ the number of maximal consecutive abstractions of length $k$, except the one at the root.

Now we introduce some statistics on degree trees. Given $(T, \ell)$ a degree tree, we denote by $\lnode(T, \ell)$ the number of leaves in $T$, by $\znode(T, \ell)$ the number of internal nodes whose leftmost descending edge has label $0$, by $\rlabel(T,\ell)$ the label of the root in $(T, \ell)$, and by $\edge_k(T, \ell)$ the number of edges with label $k \geq 1$. We then have the following transfer of statistics.

\begin{prop} \label{prop:stat-lambda-tree}
  Given $S \in \redskl_n$, let $(T, \ell) = \varphi(S)$. Then we have
  \begin{align*}
    \applv(S) = \lnode(T, \ell), &\quad \appla(S) = \znode(T, \ell), \\
    \excess(S) = \rlabel(T, \ell), &\quad \uchain_k(S) = \edge_k(T, \ell).
  \end{align*}
\end{prop}
\begin{proof}
  By Definition~\ref{defn:redskl-degree-tree}, the right child of a binary node $v$ in $S$ becomes the leftmost child of $v$ in $T$, meaning that $v$ is a leaf in $T$ if and only if its right child in $v$ is a leaf. We thus have $\applv(S) = \lnode(T, \ell)$ and $\appla(S) = \znode(T, \ell)$. The proof of Proposition~\ref{prop:bij-redskl-degree-tree} shows that $\excess(S) = \rlabel(T, \ell)$, and $\uchain_k(S) = \edge_k(T, \ell)$ is a direct consequence of Definition~\ref{defn:redskl-degree-tree}.
\end{proof}

The statistics mentioned in Proposition~\ref{prop:stat-lambda-tree} are also transferred to bipartite planar maps. More precisely, let $\bipartite_n$ be the set of bipartite planar maps with $n$ edges. For $M \in \bipartite_n$, we denote by $\white(M)$ (resp. $\black(M)$) the number of white (resp. black) vertices in $M$, and $\outdeg(M)$ the half-degree of the outer face. Moreover, for $k \geq 1$, we denote by $\face_k(M)$ the number of inner faces of degree $2k$. We have the following result from \cite{fang-bipartite}.

\begin{prop}[See Proposition~3.3,~3.6 and Remark~1 in \cite{fang-bipartite}] \label{prop:bij-tree-map}
  There is a bijection $\mu$ from $\mathcal{T}_n$ to $\bipartite_n$. Moreover, for $(T, \ell) \in \mathcal{T}_n$, let $M = \mu(T, \ell)$, and we have
  \begin{align*}
    \lnode(T, \ell) = \white(M), &\quad \znode(T, \ell) = \black(M), \\
    \rlabel(T, \ell) = \outdeg(M), &\quad \edge_k(T, \ell) = \face_k(M).
  \end{align*}
\end{prop}

Combining the bijections from $\redskl_n$ to $\mathcal{T}_{n-2}$ and from $\mathcal{T}_{n-2}$ to $\bipartite_{n-2}$, we have the bijection in Theorem~\ref{thm:main-3-conn}.

\begin{proof}[Proof of Theorem~\ref{thm:main-3-conn}]
  Consider $\mu \circ \varphi$ with $\varphi$ from Definition~\ref{defn:redskl-degree-tree} and $\mu$ from Proposition~\ref{prop:bij-tree-map}. It is a bijection from $\redskl_n$ to $\bipartite_{n-2}$ according to Proposition~\ref{prop:bij-redskl-degree-tree}~and~\ref{prop:bij-tree-map}.
\end{proof}

Combining Proposition~\ref{prop:stat-lambda-tree}~and~\ref{prop:bij-tree-map}, we have the following transfer of statistics from reduced skeletons of $\lambda$-terms in $\conn^{(3)}$ to planar bipartite maps.

\begin{thm} \label{thm:stat-lambda-map}
  Let $S \in \redskl_n$ and $M = \mu(\varphi(S)) \in \bipartite_{n-2}$. Then we have
  \begin{align*}
    \applv(S) = \white(M), &\quad \appla(S) = \black(M), \\
    \excess(S) = \outdeg(M), &\quad \uchain_k(S) = \face_k(M).
  \end{align*}
\end{thm}

Let $f_{\redskl}(t,x;p_1, p_2, \ldots)$ be the generating function of $\redskl$ with $t$ marking the size, $x$ marking the excess, and $p_k$ marking the number of maximal chains of unary nodes of length $k$. Let $f_{\bipartite}(t;p_1,p_2,\ldots)$ be the generating function of bipartite planar maps with $t$ marking the size, $x$ marking the half-degree of the outer face, and $p_k$ marking the number of inner faces of degree $2k$. We can thus express $f_{\redskl}$ with $f_{\bipartite}$.

\begin{prop} \label{prop:gf-chain}
   The generating functions $f_{\redskl}$ and $f_{\bipartite}$ are related by
  \[
    f_{\redskl}(t,x; p_1, p_2, \ldots) = t^2 f_{\bipartite}(t,x; p_1, p_2, \ldots).
  \]
\end{prop}
\begin{proof}
  This is a direct translation of Theorem~\ref{thm:stat-lambda-map} in terms of generating function.
\end{proof}

The following expression of $f_{\bipartite}$ was first found implicitly in \cite{BMS}, then written in the following form in \cite{BDFG}:
\begin{equation}
  \label{eq:gf-bipartite}
  f_{\bipartite}(t,x; p_1, p_2, \ldots) = (1+uz) \left( 1 - \sum_{k \geq 1} p_k z^k \sum_{\ell=1}^{k-1} u^\ell z^\ell \binom{2k-1}{k+\ell} \right), 
\end{equation}
where $z$ and $u$ are defined by the equations
\begin{equation}
  \label{eq:u-z-bipartite}
  z = t \left( 1 + \sum_{k \geq 1} \binom{2k-1}{k} p_k z^k \right), \quad u = x(1+uz)^2.
\end{equation}
By Proposition~\ref{prop:gf-chain}, we have the same expression for $t^{-2} f_{\redskl}(t,x;p_1,p_2,\ldots)$.

Furthermore, we can also use Theorem~\ref{thm:stat-lambda-map} to transfer known results on the distribution of various statistics on bipartite planar maps to $\lambda$-terms in $\skl^{(3)}$. As an example, in \cite{bipartite-root-degree}, Liskovets studied the asymptotic vertex degree distribution of many families of planar maps, including Eulerian planar maps, which are duals of bipartite planar maps. The following result was given in \cite{bipartite-root-degree}.

\begin{prop}[\cite{bipartite-root-degree}, Section~(1.3)~and~(3.1.1)] \label{prop:bipartite-root-degree}
  Let $X_n$ be the half-degree of the outer face of a bipartite planar map chosen uniformly from $\bipartite_n$. Then, for $k \geq 1$, when $n \to \infty$, we have
  \[
    \mathbb{P}[X_n = k] \to \frac{k}{3} \binom{2k}{k} \left( \frac{3}{16} \right)^k.
  \]
\end{prop}

Combining Proposition~\ref{prop:bipartite-root-degree} with Theorem~\ref{thm:stat-lambda-map}, we have the following corollary.

\begin{coro} \label{coro:bipartite-root-degree}
  Let $X_n$ be the number of abstractions at the beginning of a $\lambda$-term chosen uniformly from $\conn^{(3)}_n$. Then, for $k \geq 2$, when $n \to \infty$, we have
  \[
    \mathbb{P}[X_n = k] \to \frac{k-1}{3} \binom{2k-2}{k-1} \left( \frac{3}{16} \right)^{k-1}.
  \]
\end{coro}

\begin{rmk}
As both the bijection $\varphi$ from reduced skeletons to degree trees and the bijection $\mu$ (denoted by $\mathbf{T}_{\mathcal{M}}$ in \cite{fang-bipartite}) from degree trees to bipartite planar maps are direct, their composition is also a direct bijection. Moreover, we can take a shortcut in this composition by staying with the edge labeling when working on degree trees.
\end{rmk}

\section{Bijection from $\lambda$-terms in $\conn^{(1)}$ and $\conn^{(2)}$ to v-trees} \label{sec:bij-2}

In the following, we consider skeletons in $\skl^{(1)}$ and $\skl^{(2)}$. We start by defining a transformation $\psi$ from these skeletons to some trees with labels on their nodes.

\begin{defn} \label{defn:psi}
Given $S \in \skl^{(1)}$, We define $\psi(S)$ to be $(T, \ell)$ obtained as follows:
\begin{enumerate}
\item We label each binary node $u$ by $\excess(S_v)$, where $v$ is the right child of $u$, and remove all leaves and unary nodes. We obtain a binary tree $S'$ with node labeling $\ell'$.
\item We perform the classical bijection (\emph{switching left and right}) from binary tree $S'$ to plane tree $T$ as in the second step of the definition of $\varphi$ in Definition~\ref{defn:redskl-degree-tree}. We keep the labels of $\ell'$ on $S'$ in $\ell$ on $T$.
\end{enumerate}
For the reverse direction, we simply reverse all the operations above, then insert unary nodes \emph{from bottom to top} according to labels, but \emph{only to the right branch} of binary nodes, in accordance with (\textbf{Normality}).
\end{defn}

\begin{figure}
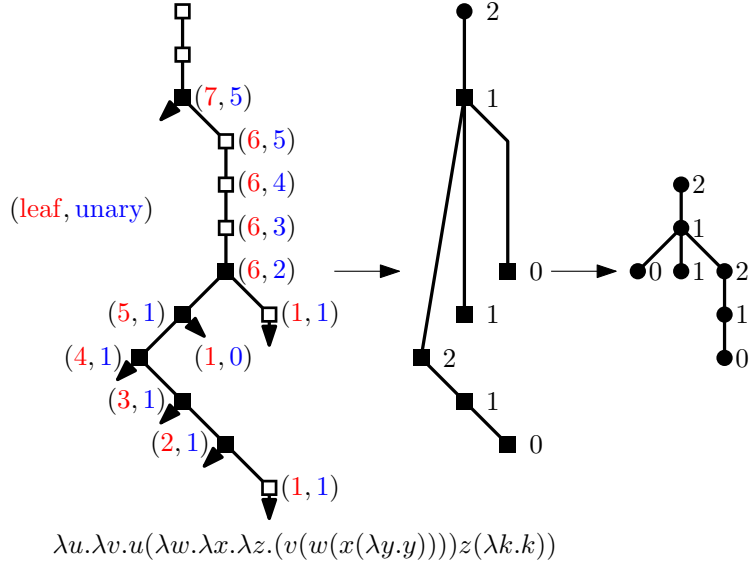

  \centering
  \insertfig{5}{0.6}
  \caption{Example of the bijection $\psi$ from skeletons in $\skl^{(1)}$ to v-trees.}
  \label{fig:v-tree-to-skeleton}
\end{figure}

We now characterize the labeled plane trees obtained by $\psi$. A node-labeled tree $(T, \ell)$ is a \tdef{v-tree} if the following conditions hold:
\begin{enumerate}[label=(\roman*)]
\item Leaves of $T$ are labeled by $0$ or $1$;\label{vtree-1}
\item Every non-root node $u$ with children $v_1, \ldots, v_k$ satisfies $0 \leq \ell(u) \leq 1 + \sum_{i=1}^k \ell(v_i)$;\label{vtree-2}
\item For the root $r$ with children $v_1, \ldots, v_k$, we have $\ell(r) = 1 + \sum_{i=1}^k \ell(v_i)$.\label{vtree-3}
\end{enumerate}
The \tdef{size} of a v-tree is its number of edges. We denote by $\vtree_n$ the set of v-trees of size $n$ and $\vtree = \bigcup_{n \geq 1} \vtree_n$. Furthermore, we denote by $\vtree^+_n$ the set of v-trees of size $n$ \emph{without label 0}, and also $\vtree^+ = \bigcup_{n \geq 1} \vtree^+_n$. We have the following characterization of $\psi(\skl^{(1)}_n)$.

\begin{prop} \label{prop:psi-chara}
  A unary-binary tree $S$ is in $\skl^{(1)}_n$ if and only if $T = \psi(S)$ is a v-tree of size $n-1$. In other words, $\psi$ is a bijection from $\skl^{(1)}_n$ to $\vtree_{n-1}$.
\end{prop}
\begin{proof}
  The correspondence in size is clear, as $S \in \skl^{(1)}_n$ has $n-1$ binary nodes, which are transformed into $n-1$ non-root nodes in $T$. In the following, by abuse of notation, we identify binary nodes in $S$ and non-root nodes in $T$.
  
  For the ``only if'' part, we first show that $T$ satisfies \ref{vtree-1} in the definition of v-trees. Let $u$ be a leaf in $T$, then it is a binary node in $S$ whose right sub-tree has no binary node, thus with only one leaf. By (\textbf{Connectedness}), this sub-tree is either a leaf or a unary node with a leaf as its child, leading to a label $0$ or $1$ of $u$ in $T$. Thus, \ref{vtree-1} is satisfied in $T$. Now we show that $T$ satisfies \ref{vtree-2}. Let $u$ be a node in $T$ that is not a leaf. With (\textbf{Connectedness}) applied to the right child of $u$, we have $\ell(u) \geq 0$. Now, in the right sub-tree of $u$ in $S$, there is at least one binary node, and we take $v_1$ to be the one closest to $u$. By (\textbf{Normality}), the leftmost path starting from $v_1$ in $S$ consists of only binary nodes $v_1, v_2, \ldots, v_k$, all transformed into children of $u$ in $T$, except a leaf $v_*$ at the end. Every unary node and leaf in the sub-tree of $S$ induced by $v_1$ is thus either $v_*$ or in one of the right sub-tree of some $v_i$. We thus have $\ell(u) \leq 1 + \sum_{i=0}^k \ell(v_i)$, with the one standing for the leaf $v_*$, and the inequality from the fact that we have not yet accounted for eventual unary nodes between $u$ and $v_1$. Therefore, $T$ also satisfies \ref{vtree-2}. $T$ also satisfies \ref{vtree-3} by (\textbf{Connectedness}) applied to the root. We thus have $T \in \vtree$.

  For the ``if'' part, we first need to show that the reverse direction of $\psi$ is always feasible on a v-tree $(T, \ell)$ to obtain a unary-binary tree $S$. The only problem that may occur is that, when we want to insert unary nodes to the right branch of a binary node $u$ with $v$ its right child, if $\ell(u) > \excess(S_v)$, then such insertion is impossible, as it may only decrease $\excess(S_v)$. However, this never happens, as from Definition~\ref{defn:psi} and the conditions of v-trees we have $\excess(S_v) = 1 + \sum_{i=1}^k \ell(v_k) \geq \ell(u)$, where $v_1, \ldots, v_k$ are the children of $u$ in $T$, with $v_1 = v$, and the $1$ stands for the leaf added to the left of $v_k$. We thus know that the reverse direction of $\psi$ from a v-tree $(T, \ell)$ is feasible. Now, the tree $S$ we obtained is indeed in $\skl^{(1)}_n$, as (\textbf{Linearity}) and (\textbf{Connectedness}) are ensured by the fact that labels in a v-tree are positive, and (\textbf{Normality}) is ensured by Definition~\ref{defn:psi}.
\end{proof}

The bijection $\psi$ specializes naturally to skeletons in $\skl^{(2)}$.

\begin{prop} \label{prop:psi-chara-2-conn}
  A unary-binary tree $S$ is in $\skl^{(2)}_n$ if and only if $T = \psi(S)$ is a v-tree of size $n-1$ without label $0$. In other words, $\psi$ is a bijection from $\skl^{(2)}_n$ to $\vtree^+_{n-1}$.
\end{prop}
\begin{proof}
  We only need to show that (\textbf{2-connectedness}) in $S$ is equivalent to the absence of label $0$ in $T$. Suppose that $u$ is a non-root non-unary node in $S$, with $k$ consecutive unary nodes above it leading to a binary node $v$. By (\textbf{Normality}), $u$ is in the right sub-tree of $v$ in $S$, meaning that $\ell(v) = \excess(S_u) - k$. We thus see that $u$ satisfies (\textbf{2-connectedness}) if and only if $\ell(u) > 0$, and we have the equivalence.
\end{proof}

\begin{rmk}
  We note that, unlike degree trees, the definition of v-tree does not distinguish left and right. Therefore, in the definition of $\psi$, there is no need to switch left and right in the second step. However, we keep the current definition for consistency with $\varphi$.
\end{rmk}

\begin{rmk}
  We notice that the definition of v-trees is very close to that of $\beta(0,1)$-trees \cite{desc-tree, description-tree}, which are in bijection with bipartite maps. The only difference is that in v-trees we may have label $1$ on leaves. However, a simple restriction forbidding label $1$ on leaves of v-trees does not lead to $\lambda$-terms in $\conn^{(3)}$, as the closed sub-term $\lambda x.x$ is translated to a leaf with label $0$ by $\psi$. Motivated by v-trees, we may also generalize $\beta(a, b)$-trees to allow extra possible labels on leaves. It is clear that the generating function of such labeled trees satisfies a functional equation with one catalytic variable, thus is algebraic \cite{BMJ}.
\end{rmk}

\section{One-corner decomposition of planar maps} \label{sec:map-decomp}

The $v$-trees also describe a recursive decomposition of planar maps, which uses the number of vertices on the outer face \emph{without multiplicity} as auxiliary statistics. This decomposition is inspired by a decomposition of bicubic planar maps \cite{desc-tree}, which is translated to a similar decomposition for bipartite maps via a bijection in \cite{Tutte-census}.

We denote by $\planar_n$ the set of planar maps with $n$ edges, and by $\planar$ the set of all planar maps. For $M$ a planar map, we denote by $\outvertex(M)$ the number of vertices of $M$ adjacent to the outer face. For the decomposition, we consider a family of planar maps called \tdef{one-corner components}, which are planar maps whose root corner is the only corner of the root vertex that is on the outer face. By definition, the empty map is also a one-corner component. We denote by $\onecorner_n$ the set of one-corner components with $n$ edges, and by $\onecorner$ the set of all one-corner components. Here, the letter $\onecorner$ stands for ``un'', which is ``one'' in French. For a one-corner component $U$, we denote by $\outvertex_{\onecorner}(U)$ the number of vertices of $U$ adjacent to the outer face, \emph{except the root vertex}. We thus have $\outvertex_{\onecorner}(U) = \outvertex(U) - 1$ for $U \in \onecorner$. Given a non-empty one-corner composition $U$ with root vertex $v$ and root edge $e$ (which may be a loop), we define $\Pi(U)$ to be the planar map obtained by deleting $e$ and re-rooting the rest at the corner from which the other end of $e$ stems.

\begin{figure}
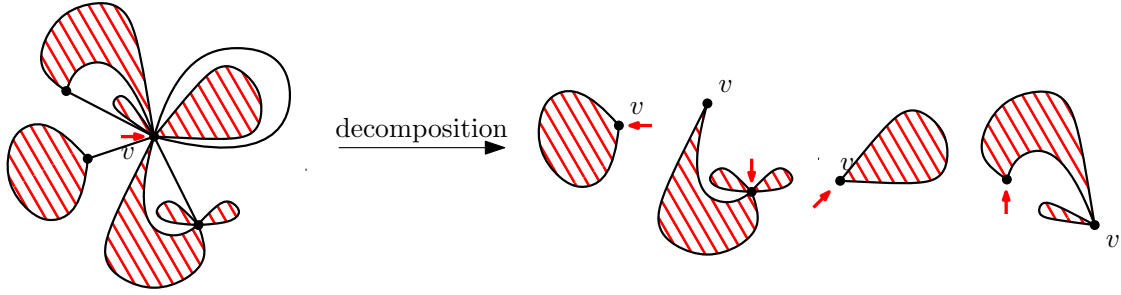

  \centering
  \insertfig{2}{0.9}
  \caption{Illustration of one-corner decomposition of planar maps}
  \label{fig:map-decomp}
\end{figure}

Given a planar map $M$, it either is empty or has at least one edge. In the latter case, we perform a decomposition as illustrated in Figure~\ref{fig:map-decomp}. Let $u$ be the root vertex of $M$. We break $M$ by cutting it at each corner of $u$ adjacent to the outer face, while duplicating $u$ for each piece. We get a sequence of one-corner component with $u$ as the root vertex, starting from the root corner in the counter-clockwise direction. Then, for each one-corner component $U_i$, we take
$M_i = \Pi(U_i)$. We thus decompose the original map into several smaller planar maps $M_1, \ldots, M_k$, and we say that it is the \tdef{one-corner decomposition} of $M$.

To reverse the decomposition, given a planar map $M$, we need to know what are the one-corner components $U$ such that $\Pi(U) = M$. We have the following result.

\begin{figure}
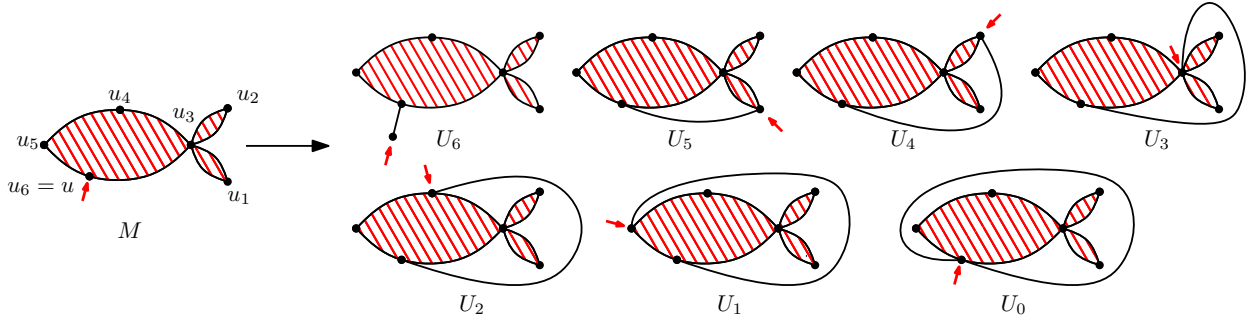

  \centering
  \insertfig{3}{1}
  \caption{Example of construction of one-corner components $U_i$ with $\Pi(U_i) = M$ for a given planar map $M$.}
  \label{fig:one-corner-recover}
\end{figure}

\begin{prop} \label{prop:one-corner-recover}
  Given a planar map $M$ with $n$ edges and $\outvertex(M) = k$, there are exactly $k + 1$ one-corner components $U_0, U_1, \ldots, U_k$, each with $n + 1$ edges, such that $\Pi(U_i) = M$ and $\outvertex_{\onecorner}(U_i) = i$ for each $0 \leq i \leq k$.
\end{prop}
\begin{proof}
  Suppose that we have $U \in \onecorner$ such that $M = \Pi(U)$. It is clear from the definition of $\Pi$ that $U$ has $n+1$ edges. The edge $e$ of $U$ that gets deleted in $M$ either is a bridge or separates the outer face with another face $f$. When $e$ is a bridge, we have $\outvertex_{\onecorner}(U) = k$, thus $U = U_k$ in this case. Otherwise, let $u_1, \ldots, u_k$ be the vertices of $M$ adjacent to the outer face, ordered by their last adjacent corner in the counter-clockwise contour of the outer face starting from the root corner. It is clear that $u_k = u$. By planarity, the other end of $e$ must be connected to some $u_i$. By the definition of $\Pi$, in order for $M = \Pi(U)$, the edge $e$ must start from the root corner of $M$ and travel counter-clockwise to a corner of $u_i$. Furthermore, the destination corner must be the last corner of $u_i$ on the outer face in counter-clockwise direction, otherwise the map obtained will not be a one-corner component. In this case, $U$ has $k - i + 1$ vertices on the outer face, namely $u_i, \ldots, u_k$. Therefore, we have $U = U_{k-i}$. These are all the possibilities. See Figure~\ref{fig:one-corner-recover} for an example of this construction.
\end{proof}

Now, to recover $M$ from its one-corner decomposition $M_1, \ldots, M_k$, by Proposition~\ref{prop:one-corner-recover}, for each $M_i$, its corresponding one-corner component $U_i$ has $\outvertex(M_i) + 1$ possibilities, with the choice given by $\outvertex_{\onecorner}(U_i)$. Gathering such statistics is sufficient to recover $M$. Using a similar approach to that in \cite{desc-tree}, we may encode a map $M$ by decomposing it recursively using the one-corner decomposition, recording the decomposition tree and the extra statistic $\outvertex_{\onecorner}(U)$ for each one-corner component $U$ appeared in the decomposition. We thus have the following theorem.

\begin{thm} \label{thm:one-corner-tree}
  For all $n \geq 0$, there is a bijection $\rho$ between the set $\planar_n$ of planar maps with $n$ edges and the set $\vtree_n$ of v-trees with $n$ edges. Furthermore, for $M \in \planar_n$, the label of the root of $\rho(M)$ is $\outvertex(M)$.
\end{thm}
\begin{proof}
  We construct $\rho$ by induction on $n$. The case $n=0$ is trivial. For the case $n=1$, we have two planar maps $M_\circ, M_-$ in $\planar_1$, where $M_\circ$ is a loop and $M_-$ is an edge between two vertices. We also have two v-trees $(T_0, \ell_0), (T_1, \ell_1)$ in $\vtree_1$, where $(T_0, \ell_0)$ is the tree with one edge with all labels $0$, while $(T_1, \ell_1)$ is the same tree with all labels $1$. We thus take $\rho(M_\circ) = (T_0, \ell_0)$ and $\rho(M_-) = (T_1, \ell_1)$, which satisfy our claim.

  Suppose that we have constructed the bijection $\rho$ for all $n' < n$. For $M \in \planar_n$, let $v_M$ be its root vertex. We split $M$ at all corners of $v_M$ adjacent to the outer face. Let $U_1, \ldots, U_k$ be the one-corner components thus obtained in counter-clockwise order starting from the root corner. We have
  \[
    \outvertex(M) = 1 + \sum_{i=1}^k \outvertex_\onecorner(U_i),
  \]
  as a vertex of $M$ adjacent to the outer face is either $v_M$ or in some $U_i$.

  We now take $M_i = \Pi(U_i)$ for $1 \leq i \leq k$, and $n_i$ the size of $M_i$. As $\Pi$ removes exactly one edge from $U_i$, we have $n_1 + \cdots + n_k = n - k$, thus $n_i \leq n-1$ for all $i$. By induction hypothesis, $\rho(M_i)$ has been constructed. Let $v_i$ be the root of $\rho(M_i)$. To construct $\rho(M) = (T, \ell)$, we put the trees $\rho(M_1), \ldots, \rho(M_k)$ from left to right, then put a new vertex $v$ as the root, making $v_1, \ldots, v_k$ children of $v$. Furthermore, we keep the labels in each $\rho(M_i)$ except for $v_i$, where we put $\ell(v_i) = \outvertex_\onecorner(U_i)$. The label of $v$ is $\outvertex(M)$. We clearly have $\ell(v) = 1 + \sum_{i=1}^k \ell(v_i)$ as shown above. Moreover, we have $0 \leq \ell(v_i) \leq \outvertex(M_i)$, where $\outvertex(M_i)$ is the root label of $v_i$ in $\rho(M_i)$. We thus check with Proposition~\ref{prop:one-corner-recover} that all the conditions of v-tree are satisfied by $\rho(M)$.

  To show that $\rho$ is indeed a bijection for the case $n$, we now show that, given a v-tree $(T, \ell) \in \vtree_n$, there is a unique $M \in \planar_n$ such that $\rho(M) = (T, \ell)$. Let $v_1, \ldots, v_k$ be the children of the root of $T$ and $T_1, \ldots, T_k$ the sub-trees of $T$ rooted at $v_1, \ldots v_k$ respectively. For each $T_i$, we construct the labeling $\ell_i$ that agrees with $\ell$ except for $v_i$, and we take $\ell_i(v_i)$ to be one plus the sum of the labels of its children. It is clear that $(T_i, \ell_i)$ is a v-tree of size strictly smaller than $n$. By induction hypothesis, there is a unique planar map $M_i$ such that $\rho(M_i) = (T_i, \ell_i)$ and $\ell_i(v_i) = \outvertex(M_i)$. We then construct the one-corner component $U_i$ with $\outvertex_\onecorner(U_i) = \ell(v_i)$. As $0 \leq \ell(v_i) \leq \ell_i(v_i) = \outvertex(M_i)$, by Proposition~\ref{prop:one-corner-recover}, such $U_i$ always exists. We then glue the marked corner of each $U_i$ to a new vertex $v_M$, from $U_1$ to $U_k$ in counter-clockwise order, to obtain a planar map $M$ with the corner between $U_1$ and $U_k$ (or the marked corner of $U_1$ when $k=1$) being the root corner. It is clear from construction that $\rho(M) = (T, \ell)$. The uniqueness of $M$ is given by the uniqueness of each $U_i$ ensured by Proposition~\ref{prop:one-corner-recover}. We thus conclude the induction.
\end{proof}

By the proof of Theorem~\ref{thm:one-corner-tree}, the v-tree $\rho(M)$ of a planar map $M$ is in fact the recursive decomposition tree of $M$ under the one-corner decomposition. The bijection $\rho$ is recursively defined in Theorem~\ref{thm:one-corner-tree}, which is less satisfying than a direct bijection. We now define a direct version of $\rho$ using an exploration procedure of planar maps, in the vein of \cite{nonsep, sticky, fang-bipartite}.

\begin{figure}
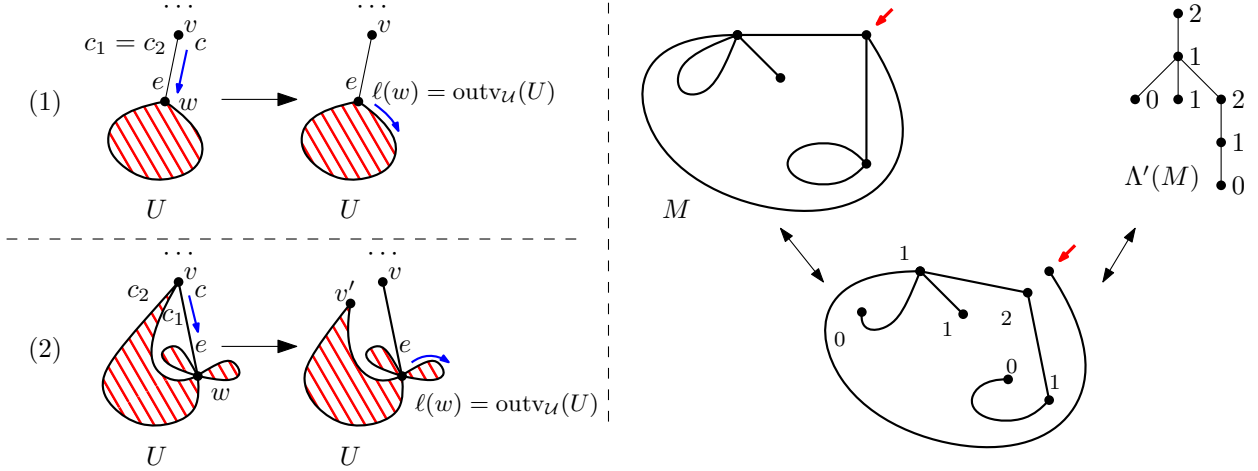

  \centering
  \insertfig{4}{1}
  \caption{Definition and example of the direct bijection $\rho'$.}
  \label{fig:map-to-v-tree-direct}
\end{figure}

Given a planar map $M$, we put the label $\outvertex(M)$ on its root vertex, and then start a contour walk in the clockwise direction from the root corner. When we walk along an edge already visited, we do nothing. Suppose that we are at a corner $c$ and is about to walk along a new edge $e$ starting from $v$ to $w$. We find the corner $c_1$ of $v$ on the other side of $e$ and the next corner $c_2$ in clockwise direction of $v$ that is adjacent to the outer face. The part of $M$ between $c$ and $c_2$ is a one-corner component $U$. We put the label $\outvertex_\onecorner(U)$ on $w$. Now we have two cases, as illustrated on the left part of Figure~\ref{fig:map-to-v-tree-direct}:
\begin{enumerate}
\item If $c_1 = c_2$, then $e$ is a bridge, and we do nothing.
\item If $c_1 \neq c_2$, then $e$ separates the outer face with an inner face. In this case, we detach the part $M'$ of $M$ between $c_1$ and $c_2$ while duplicating the vertex $v$.
\end{enumerate}
In both cases, planarity and size are conserved. After the whole walk, we obtain a planar map with labels on its vertices, denoted by $\rho'(M)$. The right part of Figure~\ref{fig:map-to-v-tree-direct} gives an example of the construction of $\rho'(M)$. It is in fact always a tree rooted at the original root corner of $M$, as every inner face will be eventually ``opened up'' in the second case above. The following theorem shows that $\rho'$ is exactly $\rho$.

\begin{thm} \label{thm:one-corner-tree-direct}
  For all $n \geq 0$ and $M \in \planar_n$, we have $\rho(M) = \rho'(M)$.
\end{thm}
\begin{proof}
  We proceed by induction on $n$. The cases $n = 0, 1$ are trivial. Suppose that our claim holds for any $n' < n$. Given $M \in \planar_n$, let $v_M$ be its root vertex, and $U_1, \ldots, U_k$ the one-corner components obtained by splitting $v_M$ by all its corners adjacent to the outer face, counter-clockwise ordered, starting from the root corner. Let $M_i = \Pi(U_i)$ for $1 \leq i \leq k$. From the construction of $\rho'(M)$, we see that the exploration starts from $U_k$, and the first visited edge $e$ becomes a bridge linking $v_M$ to the root corner of $M_k$ at some vertex $v_k$, with the size of $M_k$ strictly smaller than $n$. By induction hypothesis, the remaining part of the visit will convert $M_k$ into $\rho(M_k)$, rooted at $v_k$, with $v_k$ labeled by $\outvertex_\onecorner(U_k)$. The same applies to all components from $U_{k-1}$ to $U_1$. We thus have $\rho'(M) = (T, \ell)$, with $T$ rooted at $v_M$, with sub-trees $\rho(M_1), \ldots \rho(M_k)$ from left to right, except that the root $v_i$ of $\rho(M_i)$ is labeled by $\ell(v_i) = \outvertex_\onecorner(U_i)$. The root $v_M$ is labeled by $\ell(v_M) = \outvertex(M)$, which is equal to $1 + \sum_{i=1}^k \outvertex_\onecorner(U_i) = 1 + \sum_{i=1}^k \ell(v_i)$. By the construction of $\rho(M)$ in the proof of Theorem~\ref{thm:one-corner-tree}, we have $\rho'(M) = \rho(M)$, thus concluding the induction.
\end{proof}

\begin{proof}[Proof of Theorem~\ref{thm:main-conn}]
  As we may replace $\lambda$-terms in $\conn^{(1)}$ with their skeletons in $\skl^{(1)}$, the bijection is given by $\psi \circ \rho$. The first part is a consequence of Proposition~\ref{prop:psi-chara} and Theorem~\ref{thm:one-corner-tree}. For the second part, we notice that the existence of a loop in a planar map $M$ is equivalent to that of a one-corner component $U$ with $\outvertex_\onecorner(U) = 0$ in the recursive construction of $\rho(M)$, meaning the existence of a label $0$ in $\rho(M)$. Therefore, $\rho(\vtree^+)$ is exactly the set of loopless planar maps. We thus conclude by Proposition~\ref{prop:psi-chara-2-conn}.
\end{proof}

\begin{rmk}
  Zeilberger and Giorgetti has given a bijection between planar maps and planar linear normal $\lambda$-terms in \cite{planar-lambda}, using the Tutte decomposition \cite{Tutte-census}. Furthermore, Zeilberger also showed in \cite{lambda-trivalent} that the same bijection sends terms that are also unitless to bridgeless planar maps, which are the duals of loopless planar maps. We have checked that our bijection is different from the one in \cite{planar-lambda}, even after taking the dual. This is to be expected, as we use a different recursive decomposition. However, as there are many different conventions we may take when defining our bijections, we cannot exclude the possibility that some special choices would recover the bijection in \cite{planar-lambda}.
\end{rmk}

\begin{rmk}
  Like degree trees, we may also define an edge labeling for v-trees by putting the difference between the label of a node $u$ and its maximal possible label (1 plus the sum of the labels of its children). This presentation makes v-trees seem much closer to degree trees. As in the case of degree trees and bijectively related bipartite maps (see, \textit{e.g.}, \cite{cf-bipartite}), we may write a generating function $f_\vtree(t, x; q_1, q_2, \ldots)$ of v-trees with infinitely many variables $p_1, p_2, \ldots$, where $p_k$ marks the number of edges with label $k$. Using the standard symbolic method, we have
\[
  f_\vtree(t, x) = x + t f_\vtree(t, x) (1 + \Omega) f_\vtree(t, x),
\]
with the linear operator $\Omega$ defined by $\Omega x^n = \sum_{k=1}^n p_k x^{n-k}$ for all $n \geq 0$. We see that setting $p_k=1$ for all $k$ turns $\Omega$ into the divided difference $f(x) \mapsto \frac{f(x) - f(1)}{x - 1}$. For comparison, we have the following functional equation for the generating function of bipartite maps (see, \textit{e.g.}, \cite{cf-bipartite}), thus also of degree trees:
\[
  f_\bipartite(t, x) = 1 + t x (f_\bipartite(t, x) + \Omega) f_\bipartite(t, x).
\]
However, face degrees in bipartite maps enjoy an algebraic interpretation in the symmetric group (see, \textit{e.g.}, \cite{cf-bipartite}), which lends important structure that is absent from v-trees. Furthermore, the combinatorial interpretation of edge labels of v-trees on bipartite maps does not seem to be really natural, and the coefficients of monomials in $f_\vtree$ does not seem to be nice. These weakens further the hope of getting a nice expression for $f_\vtree(t, x)$.
\end{rmk}

\bibliographystyle{alpha}
\bibliography{lambda-bipartite}

\end{document}